\documentclass[11pt]{article}
\usepackage{latexsym,amsthm,verbatim,ifthen,graphicx,color,mathrsfs}
\usepackage[all]{xy}
\usepackage{color}

\usepackage{lipsum}

\newcommand\blfootnote[1]{%
  \begingroup
  \renewcommand\thefootnote{}\footnote{#1}%
  \addtocounter{footnote}{-1}%
  \endgroup
}

\usepackage{color}

\usepackage{lscape}
 \usepackage[ansinew]{inputenc}
 \usepackage{multicol}
 \usepackage[all]{xy}\xyoption{poly}
 \usepackage{mathrsfs}
\usepackage[psamsfonts]{eucal}
\usepackage{amssymb, amsmath}
\usepackage{geometry}
\usepackage{enumerate}
\usepackage{float}
\newtheorem*{theorem*}{Theorem}
\newtheorem{theorem}{Theorem}[section]

 \newtheorem{lemma}[theorem]{Lemma}
 
 \theoremstyle{definition}
 
 \theoremstyle{remark}

\def\quic{{\mathscr C}}

\newcommand{\bz}{\mathbb Z}
\def\bq{\mathbb{Q}}

\newcommand{\tor}{\operatorname{\text{Tor}}}

 \newcommand{\lib }{\mathbb{L}}

\newcommand{\Span}{\operatorname{{Span}}}

   \newcommand{\libc}{{\widehat\lib}}

\begin{document}

\title{The homology of the lamplighter Lie algebra}

\author{Yves F\'elix and Aniceto Murillo}

\maketitle

\begin{abstract}
We show that the associated Lie algebra of the Malcev $\bq$-completion of the lamplighter group is the pronilpotent completion of the lamplighter Lie algebra. We also prove that the homology of this completed Lie algebra is of uncountable dimension on each degree.
 \end{abstract}

Consider  the integer lamplighter group $G$ which can also be regarded as the restricted wreath product $\bz\wr\bz$ of two infinite cyclic groups and can be presented as
 $$G = \langle\,a,b\,\,\vert \, \,[a, a^{b^i}]= 1, \,\, i\in \mathbb Z\,\rangle.$$
The  explicit description of both, the pronilpotent and Malcev $\bq$-completion of $G$, is the cornerstone of the recent work of Sergei Ivanov and Roman Mikhailov on the homology of  the completion of free groups (\cite{IM}).

On the Lie side, we use the term {\em lamplighter Lie algebra} to name the  Lie analogue of $G$ (over any commutative ring $R$) introduced in  \cite{IMZ} as $L_R=R[x]\rtimes Rt$, the semidirect product of the abelian Lie algebras $R[x]$ and $R t$ with $[p,t]= x\cdot p$ for any polynomial $p$.  In \cite{IMZ} the authors prove  that the degree 2 homology $H_2(\widehat L_R; R)$ of the pronilpotent completion of $L_R$ is uncountable.  This is then used  as an important step in the proof of the non-countability dimension of  $H_2(F_R;R)$ where $F_R$ denotes the pronilpotent completion of the $R$-free Lie algebra on two generators.

 In this text we first prove:

 \vspace{3mm}\noindent {\bf Theorem 1.} {\em The pronilpotent completion $\widehat L$ of the rational  lamplighter Lie algebra is the natural Lie algebra associated to  $\widehat G$, the Malcev $\bq$-completion of the lamplighter group. In other terms, $\widehat G$ is isomorphic to the group $(\widehat L,*)$  where the multiplication $*$ is given by the Baker-Campbell-Hausdorff product.}

 \vspace{3mm} Our second result is an extension of  \cite[Theorem 2.11]{IMZ} over the rationals:

 \vspace{3mm}\noindent {\bf Theorem 2.} {\em For any $q\geq 2$, $H_q(\widehat{L})$ is uncountable. In particular, the cohomological dimension of $\widehat{L}$ is infinite.}

 \vspace{3mm} This result can be seen as a first step in the direction of a solution to the following general problem.

 \vspace{3mm}\noindent {\bf Problem}. Suppose $L=\varprojlim L/L^n$ is a pronilpotent Lie algebra, with dim$\, L=\infty$ and dim$\, L/[L,L]<\infty$. Is the cohomological dimension of $L$ always infinite?

  \vspace{3mm}

 \section{Malcev completion}

Given a group $G$  we denote as usual by $\{G^n\}_{n\ge1}$ its lower central series: $G^1= G$ and for $n\geq 2$, $G^n= [G, G^{n-1}]$.  A {\em Malcev group} is a group $G$ such that each $G^n/G^{n+1} $ is a $\mathbb Q$-vector space and $G \cong\varprojlim_n G/G^n$.
 The \emph{Malcev completion} of a group $G$ is a homomorphism of groups $ G\to \widehat G$ where $\widehat G$ is Malcev and which  induces isomorphisms $G^n/G^{n+1}\otimes \mathbb Q \stackrel{\cong}{\to}{\widehat G}^n/{\widehat G}^{n+1}$ for each $n\ge1$ (\cite[\S 7.4]{FHTII}). The  Malcev completion of $G$ satisfies $\widehat{G}= \varprojlim_n \widehat{G/G^n}$,

 On the other hand,  given a Lie algebra $L$, we also denote by $\{L^n\}_{n\ge1}$ its lower central series and  by $\widehat L= \varprojlim L/L^n$ its {\em pronilpotent completion}. A Lie algebra  $L$ is {\em pronilpotent} if the natural map $L\stackrel{\cong}{\to}\widehat L$ is an isomorphism.

Recall that, in general terms, the original Malcev correspondence (\cite{mal}) which establishes an isomorphism between the categories of nilpotent $\bq$-groups (i.e. nilpotent Malcev groups) and nilpotent Lie algebras, extends to an equivalence between the category of Malcev groups and that of pronilpotent Lie algebras (\cite[\S{A3}]{qui}). Through this equivalence, the associated Lie algebra to a given Malcev group $G$ can be described as  the limit of the  Lie algebras associated to $G/G^n$, for $n\ge1$.

 This correspondence is closely related to topology. In fact, one of the more direct constructions of the Malcev completion of a finitely generated group $G$ and its associated pronilpotent Lie algebra arises from rational homotopy theory whose bases can be found in \cite{FHT} and  \cite{FHTII}. Let
$$(\land W(n),d) \stackrel{\simeq}{\longrightarrow} A_{PL} B(G/G^{n+2})$$
 be the Sullivan minimal model of the classifying space $B(G/G^{n+2})$.
 Since $G/G^{n+2}$ is a nilpotent finitely generated group, $W(n)$ is finite dimensional and concentrated in degree 1. Denoting
$(\land W,d)= \varinjlim_n (\land W(n),d)$ it follows that  the induced morphism
$$ (\land W,d) \to A_{PL}(BG)$$
is a 1-minimal model. This means that there is a graded vector space $Z= Z^{\geq 2}$ and a differential on $\land W\otimes \land Z$ that makes $\land W\otimes \land Z$ a minimal model for $BG$ (\cite[Theorem 7.4]{FHTII}).

Now, denote by  $L= s^{-1}\mbox{Hom} (W,\mathbb Q)$ the homotopy Lie algebra of $(\land W,d)$. Recall that, since $d(W)\subset  \land^2W$, $L$ admits a Lie bracket defined by
$$<w, s[x,y]> = -<dw,sx, sy>,$$
where
$<w_1\land w_2, sx, sy> = <w_1, sx>.<w_2, sy> - <w_1, sy>.<w_2, sx>$ (\cite[\S 2.1]{FHTII}).

Finally let $I_L$ be the augmentation ideal of the universal enveloping algebra $UL$ of $L$, and let $\widehat{UL}= \varprojlim_n UL/I_L^n$ its $I_L$-adic completion. The algebra $\widehat{UL}$ is a Hopf algebra, whose primitive and group-like elements are denoted respectively by $P_L$ and $G_L$. Then (\cite{Serre}, \cite{FHTII}), $G_L$ is a group, $P_L$ is isomorphic to $\widehat L$  and the morphisms $\exp$ and $\log$ induce inverse bijections between $\widehat L$ and $G_L$. Then (\cite[Cor. 7.3]{FHTII}),  $G_L$ that is the Malcev completion of $G$  and $\widehat L$ is its associated Lie algebra.

Also arising from the Lie approach to rational homotopy theory, see \cite[\S12.1]{bfmt}, one finds an explicit description of the Malcev completion  and its associated Lie algebra of any finitely presented group  $
G=\langle \,a_1, \dots , a_p\,\, |\,\,   b_1,\dots,b_k\,\rangle
$. These are given by the quotient Lie algebra (and the corresponding quotient group with the BCH structure)
$\widehat G=\widehat{\mathbb L} (a_1, \dots , a_p)/ ( b_1,\dots,b_k),
$
where $\libc(a_1,\dots,a_n)$ is the pronilpotent completion of the free Lie algebra and each of the $b_i$'s is expressed with the BCH product.

As needed later on,  we also recall here the explicit Malcev correspondence in the  finitely generated nilpotent case.   Denote by $U(n)$ the group of $n\times n$ {\em strictly triangular matrices} with rational entries, i.e., matrices $M = (m_{ij})$ with $m_{ij}= 0$ if $i\ge j$.  In $U(n)$ we consider the Lie bracket given by commutators. On the other hand, denote by $T(n)$ the group of $n\times n$  {\em unitriangular matrices} also over $\bq$, that is, upper triangular matrices where all entries in the diagonal are 1. Then,  any finitely generated nilpotent $\bq$-group $G$ can be embedded in $T(n)$ for some $n$ and the restriction of the bijections
$$
\xymatrix{U(n)\,\, \ar@<0.98ex>[r]_(.53){\scriptscriptstyle \cong}^-{\exp} &{\,T(n)} \ar@<0.95ex>[l]^-{\log} }
$$
to $G$
determines   $L=\exp(G)$ which is the Lie algebra associated to $G$.

\section{The link between the lamplighter group and the lamplighter Lie algebra}

In what follows  $G = \langle a,b\, \vert \, [a, a^{b^i}]= 1, \,\, i\in \mathbb Z\rangle $ and $L=\bq[x]\rtimes \bq t$  will respectively denote again the lamplighter group and the (rational)  Lie algebra respectively.

\begin{lemma}
\label{lampLie} The Lie algebra $L/L^n$ is isomorphic to  the Lie algebra $E(n)$ of square $(n+1)\times (n+1)$-matrices of the form
$$
\left( \begin{array}{cccccc}
0 & b_0 & b_1 & b_2 & \dots & b_{n-1}\\
0 & 0 & a & 0 & \dots & 0\\
0 & 0 &0 & a & \dots & 0\\
\dots & \dots & \dots  & \dots  & \dots  & \dots  \\
0 & 0 & 0 & 0 &\dots & a\\
0 & 0 & 0 & 0 & \dots & 0\end{array} \right).$$
\end{lemma}

In other words, $E(n)$ is the sub Lie algebra of $U(n+1)$ generated by the    matrices $M= (m_{ij})$ with
$m_{2,3}= m_{i, i+1}$ for $i\geq 2$, and $m_{ij}= 0$ when $i\geq 2$ and $j\neq i+1$.

\begin{proof} First notice that the quotient Lie algebra $L/L^n$ is the semi direct product of abelian Lie algebra
$$L/L^n = \mathbb Q[x]/x^n \rtimes \mathbb Q t$$
with $[p(x), t]= x.p(x)$.

Denote by $A$  the matrix of $E(n)$ corresponding to the values $a=1$ and $b_i=0$ for all $i=0,\dots, n-1$. In the same way, for $r=0,\dots,n-1$, let $B_r\in E(n)$ corresponding to $a=0$ and $b_i=\delta_i^r$. A direct calculation shows that $AB_r = 0$ for all $r$, $B_rB_s = 0$ for all $r,s$, and  $[B_r, A]= B_rA= B_{r+1}$ for $0\le r\le n-2$ while $[B_{n-1},A]=0$.
It follows then that  the map
$$
\varphi\colon L/L^n\stackrel{\cong}{\longrightarrow} E(n),\quad\varphi (t) = A,\quad \varphi (x^r) = B_r
$$ is an isomorphism of Lie algebras.
\end{proof}

Denote by  $G(n)=\exp\bigl(E(n)\bigr)$  the subgroup of $T(n+1)$ which, by the previous Lemma,  is   the Malcev group associated to $L/L^n$.

\begin{lemma}
\label{lampgroup} $G(n)\cong\widehat{G/G^n}$.
\end{lemma}

\begin{proof} Consider the map
$$ \psi\colon G/G^n \longrightarrow G(n),\quad \psi (b) = e^{-A},\quad \psi (a)= e^{B_0} = 1+B_0.$$
From the formulas
$$\psi (a^{b^i}) = e^{-iA}e^{B_0} e^{iA}= 1+ B_0e^{iA},\quad \psi (a\cdot a^{b^i})= \psi (a^{b^i}\cdot a) = 1+B+Be^{iA},\quad i\in\bz,$$
it follows that $\psi$ is a morphism of groups.
Let $H_1$ be the completion of the abelian group generated by $a,  a^b, \dots, a^{b^{n-1}}$, and $H_2$ the vector space generated by the elements $e^{B_i}$, $i=0,\dots,n-1$. Then, we have a commutative diagram of short exact sequences:
$$\xymatrix{ 0 \ar[r] & H_1\ar[d]^\psi \ar[rr] && \widehat{G/G^n} \ar[d]^{\widehat\psi} \ar[rr] && \mathbb Q b\ar[d]^\psi \ar[r] & 0\\
0 \ar[r] &H_2\ar[rr] &&G(n)\ar[rr] && \{e^{\lambda A},\, \lambda\in \mathbb Q\} \ar[r] & 0}$$
Therefore, $\widehat\psi$ is an isomorphism of groups. \end{proof}

\begin{proof}[Proof of Theorem 1]

From Lemmas  \ref{lampLie} and  \ref{lampgroup} we get:
$$\exp (\widehat{L} ) = \varprojlim_n \exp ({L/L^n}) \cong \varprojlim_n G(n) \cong \varprojlim_n (\widehat{G/G^n}) = \widehat{G}.$$
\end{proof}

\section{Homology of the completion of the lamplighter Lie algebra}

Recall that the homology of a Lie algebra $E$ is defined as $H(E)=\tor_{UE}(\bq,\bq)$ and coincides with the homology of the {\em chain coalgebra} $\quic(E)$, see for instance \cite[p. 301]{FHT}. This  is the differential graded coalgebra $\quic(E)=(\land sE,d)$ where $s$ denotes suspension, i.e., $sE$ is concentrated in degree $1$. The differential, which vanishes on $\quic_0(E)=\bq$ and  $\quic_1(E)=sE$, is given in $\quic_p(E)$ by,
$$d (sx_1\land \dots \land sx_p)= \sum_{1\leq i<j\leq p} (-1)^{i+j}s[x_i, x_j]\land sx_1\land \dots \widehat{sx_i}\dots \widehat{sx_i}\dots \land sx_p,$$
for any $p\ge 2$.

In our situation, a generic element of the Lamplighter Lie algebra
 $L$ (resp. its pronilpotent completion $\widehat{L}$)  is a sum (resp. series) $\lambda t+\sum_{r\ge0}\lambda_r\,x^r$ with $\lambda,\lambda_r\in\bq$.

We compute first $H_2(L)$ and $H_2(\widehat{L})$. For it, note that
a basis of cycles in $\quic_2(L)$ is given by the $\{sx^r\land sx^s\}_{0\le r<s}$.  Denote
 $$
 V_p= \{ sx^r\land sx^s,\,\,\,\, 0\leq r<s,\, r+s= p\}$$
and observe that $\dim V_{2n}= n$ and $\dim V_{2n+1}= n+1$.
Moreover, the differential on $\quic_3(L)$ satisfies
 $$d (sx^r\land sx^s\land st) =sx^{r+1}\land sx^s + sx^r\land sx^{s+1},$$
 so that it restricts to morphisms
 $$d\colon  V_p\land st \to V_{p+1}.$$

\begin{lemma}
\label{calcul2} With the above notations,
\begin{enumerate}
\item[(i)] The differential $d \colon V_n\land st\to V_{n+1}$ is injective  for any $n\ge1$ and an isomorphism for $n$ odd.
\item[(ii)] A basis of $H_2(L)$ is given by the class of the cycles $sx^0\land sx^p$, with $p$ odd.
\item[(iii)] $H_2(\widehat{L})$ is isomorphic to $sx^0\land s(x\cdot\bq\bigl[[x]\bigr])$. In particular, it is of uncountable dimension.
\end{enumerate}
\end{lemma}

\begin{proof} A direct computation shows that $d \colon V_n\land st\to V_{n+1}$ is injective  for all $n$. On the other hand, note that
each cycle of $\quic_2(L)$ can be written as $sx^0\land sx^r + d(\gamma)$ for some $r$ and some $\gamma\in\quic_3(L)$. Furthermore,
 $$sx^0\land sx^{2n}= d\, \left( \sum_{i=0}^{n-1} (-1)^i \, sx^i \land sx^{2n-i-1}\land st\, \right)$$
which proves that $d\colon V_{2n-1}\land st \to V_{2n}$ is an isomorphism. As a consequence,
  a basis for $H_2(L)$ is given by the elements $sx^0\land sx^p$, with $p$ odd. This also, implies that $H_2(\widehat{L}) = sx^0 \land S$ where $S$ is the vector space of series $\sum_i \lambda_i \, sx^{2i+1}$, with $\lambda_i\in \mathbb Q$.
\end{proof}

We now analyze $H_q(\widehat L)$ for any  $q\ge 2$. The vector space $\quic_q(L)$ decomposes as a direct sum $\quic_q(L) = A_q\oplus B_q$, where
$A_q $ is the sub vector space generated by the elements $sx^{i_1}\land \dots \land sx^{i_q}$,  with $i_1<\dots <i_q,$
and
 $B_q $ is the sub vector space generated by the elements $sx^{i_1}\land \dots \land sx^{i_{q-1}}\land st$, with $ i_1<\dots <i_{q-1}.$

The differential $d \colon  \quic_{q+1}(L)\to \quic_q(L)$ is zero on $A_{q+1}$ and maps $B_{q+1}$ to $A_q$,
$$d(sx^{j_1}\land \dots \land sx^{j_q}\land st ) = \sum_{p=1}^q \pm sx^{j_1}\land \dots \land sx^{j_p+1}\land \dots \land sx^{j_q}.$$
Since $B_{q+1}= A_q\land st$, we may consider the bijection $A_q\stackrel{\cong}{\to} B_{q+1}$, $\Phi\mapsto \Phi\land st$, and will identify $d\colon B_{q+1}\to A_q$ with $A_q\stackrel{\cong}{\to} B_{q+1} \stackrel{d}{\to}A_q$.

Denote by $V_{q,n}\subset A_q$  the sub vector space generated by the elements $sx^{i_1}\land \dots \land sx^{i_q}$ with $\sum_{j=1}^q i_j = n$. For sake of simplicity we replace the elements $sx^{i_1}\land \dots sx^{i_q}$ by the sequence $(i_1, \dots i_q)$ and thus,
$$
V_{q,n}=\Span\{\,(i_1, \dots , i_q),\,\,0\leq i_1<i_2 \dots <i_q,\,\,\sum_{j=1}^q i_j = n\,\}.
$$

We decompose $V_{q,n}$ as the direct sum
$$V_{q,n} = W_{q,n}\oplus E_{q,n},$$
with the property that $(i_1, \dots , i_q)$ belong to $W_{q,n}$ when $i_1\geq 1$, and belongs to $E_{q,n}$ if $i_1= 0$. Finally let $\rho\colon V_{q,n}\to E_{q,n}$ the projection with kernel $W_{q,n}$.
With this notation we have:

\begin{lemma}
\label{calcul}
\begin{enumerate}
\item[(i)] For any $q\ge 2$ and any $n\ge1$, the linear map  $d\colon W_{q,n}\to W_{q, n+1}$ is injective.
\item[(ii)] For any $q\ge 2$ and any $n\ge1$, the linear map $\rho\circ d \colon E_{q,n}\to E_{q, n+1}$ is injective.
\end{enumerate}
\end{lemma}

\begin{proof} We proceed by induction on $q$ and for fixed $q$ by induction on $n$. We can thus suppose that the result is true for $p<q$ and when $p=q$, for $(q,n')$ with $n'<n$. Combining (i) and (ii) we deduce that $d \colon  V(p,n')\to V(p, n'+1)$ is injective for $p<q$ and when $p=q$ for $n'<n$.

The initial step of the induction for $q= 2$ is (i) of  Lemma \ref{calcul2}. Moreover, note that,
$$
V_{q,n}= 0\quad\text{for}\quad q\ge2 \quad\text{and}\quad  n<\frac{q(q-1)}{2}.
$$
 This provides (i) and (ii) for small values of $n$ and the initial induction step is set.

\medskip
Consider the isomorphism
$$\varphi \colon  W_{q,n}\stackrel{\cong}{\longrightarrow} V_{q,n-q}, \quad \varphi (i_1, \dots , i_q)= (i_1-1, \dots , i_q-1),$$
and the commutative diagram
$$\xymatrix{
W_{q,n} \ar[d]^d\ar[rr]_\cong^\varphi && V_{q, n-q}\ar[d]^d\\
W_{q, n+1}\ar[rr]_\cong^\varphi && V_{q, n+1-q}.
}$$
By induction,  $d \colon  W_{q,n}\to W_{q, n+1}$ is injective and (i) is proved.

\medskip
In the same way, (ii)  is a direct consequence of the isomorphism
$$
\psi \colon  E_{q,n}\stackrel{\cong}{\longrightarrow} V_{q-1, n-q+1},\quad\psi (0, i_1, \dots , i_q) = (i_2-1, \dots i_q-1),
$$
and the commutative diagram
\begin{equation}\label{star}\xymatrix{
E_{q,n}\ar[d]^{\rho\circ d} \ar[rr]_(.45)\cong^(.45)\psi && V_{q-1, n-q+1}\ar[d]^d\\
E_{q,n+1} \ar[rr]_(.45)\cong^(.45)\psi && V_{q-1, n-q+2}.
}
\end{equation}

\end{proof}

\begin{lemma}
\begin{enumerate}
\item[(1)] The differential $d\colon  V_{q,n}\to V_{q, n+1}$ is an injective map.
\item[(2)] When $q\geq 2$, the elements $(0, 1,2, \dots , q-2, r)$ with $q+r$ odd and $r>q-2$, are not in the image of $d$.
\end{enumerate}
\end{lemma}

\begin{proof} (1) Write any given $\alpha\in V_{q,n}$ as $\alpha = \beta +\gamma$, with $\beta\in W_{q,n}$ and $\gamma\in E_{q,n}$. If $0= d\alpha$, then $0= \rho(d\alpha)= \rho(d \gamma)$. Thus, by (ii) of  Lemma \ref{calcul}, $\gamma = 0$ and hence, by  (i) of Lemma \ref{calcul}, $\beta$ also vanishes.

\vspace{1mm} (2) We proceed by induction on $q\geq 2$. For $q= 2$ this is (ii) of  Lemma  \ref{calcul2}. Assume $q>2$ and write $\alpha = (0,1,2, \dots , q-2,r)$ with $r+q$ odd. Using the morphism $\psi$ defined above, we have
$$\psi(\alpha) = (0, 1,2, \dots , q-3, r-1).$$ By the inductive hypothesis, $\psi (\alpha)$ is not in the image of $d$ so that, by the commutativity of the diagram (\ref{star}), $\alpha$ is not in the image of $\rho\circ d\colon  E_{q,*}\to E_{q, *+1}$. Now suppose $\alpha = d(\alpha'+\beta')$ with $\alpha'\in E_{q,*}$ and $\beta'\in W_{q,*}$. Then, $\alpha =\rho( d \alpha')$  which is impossible.
\end{proof}

\vspace{3mm}Since $q^2+1 > \frac{(q+2)(q+1)}{2}$, we directly deduce  the following which in particular proves Theorem 2:

\begin{theorem} For any $q\geq 2$,
$$sx^0\land \dots sx^{q-2} \land s\bigl(x^{q^2+1}\cdot \mathbb Q\bigl[[x^2]\bigr]\bigr) $$
is a vector space of cycles that inject into $H_q(\widehat{L})$.
In particular, $H_q(\widehat{L})$ is uncountable for any $q\geq 2$.\hfill$\square$
\end{theorem}

 \end{document}